# Finite arithmetic subgroups of $GL_n$

Marcin Mazur

*Department of Mathematics, the University of Chicago,*

*5734 S. University Avenue, Chicago, IL 60637*

e-mail: mazur@math.uchicago.edu

**Abstract**

We discuss the following conjecture of Kitaoka: *if a finite subgroup $G$ of $GL_n(O_K)$ is invariant under the action of $Gal(K/\mathbb{Q})$ then it is contained in $GL_n(K^{ab})$*. Here $O_K$ is the ring of integers in a finite, Galois extension $K$ of $\mathbb{Q}$ and $K^{ab}$ is the maximal, abelian subextension of $K$. Our main result reduces this conjecture to a special case of elementary abelian $p-$groups $G$. Also, we construct some new examples which negatively answer a question of Kitaoka.

## 1  Introduction.

About 20 years ago Yoshiyuki Kitaoka, motivated by some questions concerning quadratic forms, started the investigation of finite, stable under the action of $Gal(K/\mathbb{Q})$ subgroups of $GL_n(O_K)$, where $O_K$ is the ring of integers in a finite, Galois extension $K$ of $\mathbb{Q}$. The subject seems to be important and interesting though we feel that it did not get sufficient attention. At present we have a very clear conjectural picture due to Kitaoka but the actual results are very fragmentary and the supporting evidence is rather poor. In the present article we recall the conjectures and reduce them to a very special case. Also, we describe some new examples which negatively answer one of Kitaoka's questions.

For a finite, Galois extension $K$ of the rationals let $\Gamma = Gal(K/\mathbb{Q})$ and $O_K$ be the ring of integers in $K$. In the case of totally real $K$ the main conjecture is particularly easy to state:

**Conjecture 1** *Let $G$ be a finite subgroup of $GL_n(O_K)$ stable under the action of $\Gamma$. If $K$ is totally real then $G \subseteq GL_n(\mathbb{Z})$.*

Conjecture 1 is equivalent to the following conjecture, important for the theory of arithmetic groups:

**Conjecture 2** *Let $H$ be an algebraic group of compact type defined over $\mathbb{Q}$. For any totally real number field $K$ the equality $H_\mathbb{Z} = H_{O_K}$ holds.*

In fact it would be enough to prove Conjecture 2 for orthogonal groups of positively defined quadratic forms over $\mathbb{Q}$. As a corollary we would get the following, interesting result:

*Any two positively defined quadratic forms over $\mathbb{Z}$, which are equivalent over a ring of integers of some totally real number field are already equivalent over $\mathbb{Z}$.*

For more details about the interrelation between the above conjectures and their consequences we refer to the beautiful book [6] (paragraph 4.8).

After stating Conjecture 1 it is natural to ask what happens for arbitrary Galois extensions of $\mathbb{Q}$. To spell out a conjectural answer, due to Kitaoka, consider a free $\mathbb{Z}$-module $M$ of rank $n$ with basis $e_1, ..., e_n$. The group $GL_n(O_K)$ acts in a natural way on $O_K \otimes M = \bigoplus_{i=1}^n O_K e_i$: $(a_{i,j})e_i = \sum a_{i,j}e_j$. We say that a finite subgroup $G$ of $GL_n(O_K)$ if of *A-type* if there exists a decomposition $M = \bigoplus_{i=1}^k M_i$ such that for every $g \in G$ there are a permutation $\pi(g)$ of $\{1, ..., k\}$ and roots of unity $\epsilon_i(g) \in K$ such that $\epsilon_i(g)gM_i = M_{\pi(g)}$ for $i = 1, ..., k$. It is easy to see that $\pi$ is a group homomorphism $G \longrightarrow \Sigma_k$.

**Conjecture 3** *Any finite subgroup of $GL_n(O_K)$ stable under the action of $\Gamma$ is of A-type.*

Of course for totally real $K$ Conjecture 3 reduces to Conjecture 1. The best result up to date toward Conjecture 3 is the following theorem due to Kitaoka and Suzuki ([4]):

**Theorem 1** *If the Galois group $\Gamma$ is nilpotent then Conjecture 3 holds.*

In the light of this result Conjecture 3 is easily seen to be equivalent to **Conjecture 4** *If a finite group $G \subseteq GL_n(O_K)$ is $\Gamma-$stable then in fact $G \subseteq GL_n(O_{K^{ab}})$, where $K^{ab}$ is the maximal abelian subextension of $K$.*

Our main result reduces Conjecture 4 to the case of elementary abelian $p-$groups $G$.

## 2 Diagonal $p-$groups

Suppose that $P$ is an abelian $p-$group sitting in $GL_n(\overline{\mathbb{Q}})$ as a subgroup of diagonal matrices. In other words, there are abelian characters $\chi_1, ..., \chi_n$ of $G$ such that $g = diag(\chi_1(g), ..., \chi_n(g))$. After conjugating by a permutation matrix (which has entries in $\mathbb{Z}$) we can and will assume that there are integers $k_0 = 0 < k_1 < ... < k_s < k_{s+1} = n$ such that $\chi_i = \chi_j$ iff $k_t < i, j \leq k_{t+1}$ for some $0 \leq t \leq s$. We say in this case that $P$ is in a *strongly diagonal* form. Set $N(P)$ for the normalizer of $P$ in $GL_n(\overline{\mathbb{Q}})$ and $C(P)$ for its centralizer. The following lemma is almost obvious and we omit a proof.

**Lemma 1** *The centralizer $C(P)$ of $P$ in $GL_n(\overline{\mathbb{Q}})$ equals $GL_{k_1-k_0}(\overline{\mathbb{Q}}) \times \cdots \times GL_{k_{s+1}-k_s}(\overline{\mathbb{Q}})$. The group $W(P) = N(P)/C(P)$ is finite.*

Note that $P$ is stable under the action of the absolute Galois group $\Gamma_{\mathbb{Q}}$ of $\mathbb{Q}$. Thus both $N(P)$ and $C(P)$ are $\Gamma_{\mathbb{Q}}-$stable. We denote by $\Pi_n$ the group of permutation matrices in $GL_n$. If $N \in N(P)$ and $\phi$ is the automorphism of $P$ induced by conjugation with $N$ then clearly the map $\chi \mapsto \chi \circ \phi$ permutes the characters $\chi_1, ..., \chi_n$. It is fairly obvious that we can choose a permutation $\pi \in \Sigma_n$ such that $\chi_i \circ \phi = \chi_{\pi(i)}$ and whenever $i < j$ and $\chi_{\pi(i)} = \chi_{\pi(j)}$ then $\pi(i) < \pi(j)$. Moreover such permutation is unique and if we denote by $T_N$ the corresponding permutation matrix then $NT_N^{-1} \in C(P)$ and the association $N \mapsto T_N$ is a group homomorphism $\sigma : N(P) \longrightarrow \Pi_n \cap N(P)$. We denote by $\Pi_P$ the image of $\sigma$. Thus we get the following

**Lemma 2** *The group $N(P)$ is a semidirect product of $C(P)$ and $\Pi_P$. In particular, the induced action of $\Gamma_{\mathbb{Q}}$ on $W(P)$ is trivial.*

Note that $N(P)$ acts by conjugation on $M_{k_1-k_0}(\overline{\mathbb{Q}}) \times \cdots \times M_{k_{s+1}-k_s}(\overline{\mathbb{Q}})$ which has a basis consisting of matrices with one entry equal to 1 and all others being 0. Withe respect to this basis the action of $N(P)$ defines a representation $\rho$ of $N(P)$ in $GL_{m_1^2+\ldots+m_{s+1}^2}$, where $m_i = k_i - k_{i-1}$. It is clear that this representation respects the action of $\Gamma_{\mathbb{Q}}$. Also, if both $a$ and $a^{-1}$ have entries in $O_L$ for some number field $L$, then so do $\rho(a)$ and $\rho(a^{-1})$. The kernel of $\rho$ equals to the product of centers of $M_{k_i-k_{i-1}}(\overline{\mathbb{Q}})$. In particular, we get the following

**Lemma 3** *Suppose that $G$ is a finite, $\Gamma_{\mathbb{Q}}$-invariant subgroup of $N(P) \cap GL_n(O_K)$ for some Galois extension $K$ of $\mathbb{Q}$. Then $\rho(G)$ is a finite, $\Gamma_{\mathbb{Q}}$-invariant subgroup of $GL_{m_1^2+\ldots+m_{s+1}^2}(O_K)$ and the map $\rho: G \longrightarrow \rho(G)$ commutes with the action of $\Gamma_{\mathbb{Q}}$.*

The last lemma will allow us to perform a "dirty trick", i.e. apply Conjecture 3 to a quotient of $G$ at the expense of raising $n$.

Suppose now that $K$ is a number field and that $H \subseteq GL_n(\overline{\mathbb{Q}})$ is an elementary abelian $p$-group stable under the action of the absolute Galois group $\Gamma_K$ of $K$. There exists a matrix $E \in GL_n(\overline{\mathbb{Q}})$ such that $G = E^{-1}HE$ is in a strongly diagonal form. Since $H$ is $\Gamma_K$-stable, we have $E^{-1}E^\tau \in N(G)$ for all $\tau \in \Gamma_K$. In other words, the function $\tau^{-1} \mapsto E^{-1}E^\tau$ is a 1-cocycle of $\Gamma_K$ with values in $N(G)$. Conversely, suppose that $f: \Gamma_K \longrightarrow N(G)$ is a one cocycle. Since $H^1(\Gamma_K, GL_n) = 1$, there exists a matrix $E \in GL_n(\overline{\mathbb{Q}})$ such that $f(\tau^{-1}) = E^{-1}E^\tau$ and the group $EGE^{-1}$ is clearly $\Gamma_K$-stable. We say that $L$ is *the field of definition* of $EGE^{-1}$ if it is the smallest Galois extension of $K$ such that $EGE^{-1} \subseteq GL_n(L)$. Note that if $F$ is another matrix representing $f$ then $EF^{-1} \in GL_n(K)$, so that $EGE^{-1}$ and $FGF^{-1}$ are conjugate over $K$. Clearly homologically equivalent cocycles lead to the same $GL_n(K)$-conjugacy classes of finite subgroups. Moreover, the field of definition for $EGE^{-1}$ depends only on the homomorphism $\overline{f}: \Gamma_K \longrightarrow W(G)$ obtained by composing $f$ with the natural projection of $N(G)$ onto $W(G)$. In fact, for $\tau \in \Gamma_K$ and $g \in G$ we have $(EgE^{-1})^\tau = E(E^{-1}E^\tau)g^\tau(E^{-1}E^\tau)^{-1}E^{-1}$ and the action of $E^{-1}E^\tau = f(\tau) \in N(G)$ clearly depends only on its image in $N(G)/C(G)$. Thus for the determination of possible fields of definition we can restrict our attention to homomorphisms $\Gamma_K \longrightarrow \Pi_G$. But when we would like to study the fields of definition of subgroups in $GL_n(O_{\overline{\mathbb{Q}}})$ ($O_{\overline{\mathbb{Q}}}$ is the ring of algebraic integers) the problem is much more subtle. It can a priori happen that some cocycle $f$ gives rise to $E$ such that $EGE^{-1}$ is in $GL_n(O_{\overline{\mathbb{Q}}})$ but the corresponding $\overline{f}$ does not have this property. We can not therefore restrict only to $\text{Hom}(\Gamma_K, \Pi_G)$. It would be very interesting to get a nice description of those homomorphisms $\overline{f}$ which come from a cocycle giving rise to a subgroup consisting of matrices with entries in algebraic integers. It is easy to see that Conjecture 4 for elementary abelian $p$−groups $G$ is equivalent to saying that any such homomorphism (when $K = \mathbb{Q}$) is trivial on the commutator subgroup of $\Gamma_{\mathbb{Q}}$. We will see in the next section that Conjecture 4 reduces to this special case.

We will end this section by discussing some examples. From the very early days of the conjectures there was a desire to extend them to the relative case. In other words, we would like to know which fields can be the fields of definition for finite, $\Gamma_K$−stable subgroups of $GL_n(\overline{\mathbb{Q}})$ with entries in algebraic integers. On page 260 of [3] (see also [6], paragraph 4.8) a totally real number field $K$ is constructed for which there is a finite group with the field of definition being a totally real, unramified non-trivial extension with cyclic Galois group. There was a hope for a while that any such group has its field of definition contained in the cyclotomic extension of $K$ but counterexamples have been produced in [4]. In the same paper the authors rise the question of whether the field of definition is always an abelian extension of $K$. We show that the answer is negative and the Galois group of the field of definition can be a nonabelian simple group. The idea of our construction is very simple. Consider a finite extension $L$ of $K$ of degree $n$ and let $\tau_1, ..., \tau_n$ be the embeddings $L \hookrightarrow \overline{\mathbb{Q}}$ over $K$. Let $u_1, ..., u_n$ be a basis of $L/K$. Define the matrix $U = (u_{i,j})$ by $u_{i,j} = u_i^{\tau_j}$. Plainly $U$ is nonsingular and for any $\tau \in \Gamma_K$ we have $U^\tau U^{-1} = A_\tau \in \Pi_n$. This defines a homomorphism $\psi : \Gamma_K \longrightarrow \Pi_n$ whose kernel is $\Gamma_M$, where $M$ is the Galois closure of $L/K$. Fix a prime $p$ and let $i : \Gamma_K \longrightarrow \mathbb{Z}/(p-1)\mathbb{Z}$ be the cyclotomic character which is given by the action of $\Gamma_K$ on $p$−th roots of 1. Denote by $P$ the maximal, elementary abelian $p$−subgroup of the diagonal matrices in $GL_n(\overline{K})$. The group $P_U = U^{-1}PU$ is $\Gamma_K$−stable and for $\tau \in \Gamma_K$ and $w \in P$ we have $(U^{-1}wU)^\tau = U^{-1}\psi(\tau)^{-1}w^{i(\tau)}\psi(\tau)U$. It is easy to see that $\tau$ acts trivially on $P_U$ iff both $\psi(\tau)$ and $i(\tau)$ are trivial. In other words, the field of definition for $P_U$ equals $N = M(\xi_p)$, where $\xi_p$ is a primitive $p$−th root of 1.

The only problem with the above construction is that the group $P_U$ does not consists of matrices with integral entries in general. Suppose however that $O_L$ is a free $O_K$−module and that $L/K$ is unramified. If we take for the $u_i$'s a basis for $O_L$ over $O_K$ then $U \in GL_n(O_M)$ and consequently $P_U \subseteq GL_n(O_N)$. For an explicit example we can use the results of [7]. Let $K = \mathbb{Q}(\sqrt{36497})$. This field has class number 1 and an everywhere unramified, Galois, totally real extension $L$ with Galois group $A_5$ (the alternating group of degree 5). In particular, $O_L$ is a free $O_K$−module. Our construction in this case gives a finite, $\Gamma_K$−stable subgroup of $GL_{60}(O_{L(\xi_p)})$ with field of definition $L(\xi_p)$, where $p$ is arbitrary prime number. Extending our base field to $K(\xi_p)$ we get an example with Galois group of the field of definition equal to $A_5$.

## 3 Reduction

We start by recalling Minkowski's Lemma. Let $K$ be an algebraic number field. For an ideal $I$ of the ring of integers $O_K$ we denote by $R_I$ the natural map $GL_n(O_K) \longrightarrow GL_n(O_K/I)$ and by $GL_n(O_K, I)$ its kernel. Let $\beta$ be a prime ideal of $O_K$ lying over the rational prime $p$. By $f$ we denote the ramification index of $\beta$.

**Proposition 1** *Any torsion element of $GL_n(O_K, \beta)$ has $p$−power order. If $Gl_n(O_K, \beta^k)$ contains an element of order $p^s$ then $k \leq fp^{1-s}/(p-1)$.*

This result is classical but we include a proof for the readers convenience.
*Proof:* Let $x \in GL_n(O_K, \beta)$ be an element of prime order $q$. Let $t$ be maximal such that $x \in GL_n(O_K, \beta^t)$. Writing $x = 1 + y$ we get $1 = (1+y)^q$ and equivalently

$$-qy = \sum_{i=2}^{q} \binom{q}{i} y^i$$

If $q \neq p$ then the left hand side has entries in $\beta^t$ but not all of them are in $\beta^{t+1}$. On the other hand, the right hand side has entries in $\beta^{2t}$ which leads to a contradiction. Thus $q = p$ and the left hand side of our equality has entries in $\beta^{f+t}$ but not all of them are in $\beta^{f+t+1}$. The right hand side has entries in $\beta^{\min(f+2t,pt)}$. Thus $f + t \geq \min(f + 2t, pt)$ and therefore $f + t \geq pt$, i.e. $t \leq f/(p-1)$. This proves the lemma for $s = 1$. Also, since $k \leq t$, we get that $f + k \geq pk$. This implies that

$GL_n(O_K, \beta^{pk})$ contains an element of order $p^{s-1}$ and the result follows by induction. □

**Corollary 1** *If $p \neq 2$ then $GL_n(\mathbb{Z}, p)$ is torsionfree. For $p = 2$ the group $GL_n(\mathbb{Z}, 4)$ is torsionfree and torsion elements of $GL_n(\mathbb{Z}, 2)$ are of order at most 2.*

Suppose now that $K$ is Galois over $\mathbb{Q}$ with Galois group $\Gamma$. Let $G \subseteq GL_n(O_K)$ be a finite, $\Gamma$-stable group. For a prime $\beta$ of $O_K$ let $G(\beta)$ be the kernel of $R_\beta$ restricted to $G$. We will need the following lemma due to Kitaoka and Suzuki:

**Lemma 4** *If $K$ is an abelian extension of $\mathbb{Q}$ then $G(\beta)$ can be conjugated to diagonal matrices by an element of $GL_n(\mathbb{Z})$.*

For a proof see [4], Lemma 1.

Suppose that Conjecture 4 is false and let $G$ be a counterexample of minimal order. We also assume that $n$ is minimal possible for $G$ and $K$ is minimal. In particular, no element of $\Gamma$ acts trivially on $G$. By Lemma 3 of [2] we may assume that there is a prime $p$ such that $K$ is unramified outside $p$. We denote by $I_\beta$ and $D_\beta$ the inertia and decomposition groups of a prime ideal $\beta$ of $O_K$ respectively. We have the following easy lemma:

**Lemma 5** *The Galois group $\Gamma$ is generated by the inertia subgroups $I_\beta$ for primes $\beta$ lying over $p$.*

*Proof:* Any subfield of $K$ fixed by all the inertia of primes above $p$ is unramified over $\mathbb{Q}$ hence coincides with $\mathbb{Q}$. □

Clearly $G(\beta)$ is not trivial for any prime $\beta$ over $p$. Otherwise $I_\beta$ would act trivially on $G$. Since the inertia subgroups of different primes over $p$ are conjugate, this contradicts Lemma 5.

**Lemma 6** *$G$ is generated by the subgroups $G(\beta)$, where $\beta$ runs over all primes above $p$.*

*Proof:* By Proposition 1 we get that $G(\gamma)$ is a normal $p$-subgroup of $G$ for any prime ideal $\gamma$ of $K$ lying over $p$. Thus the composition $H$ of all $G(\gamma)$, where $\gamma$ runs over the primes above $p$, is a normal, $\Gamma$-stable $p$-subgroup of $G$. Suppose that $H \neq G$. Then $H \subseteq GL_n(K^{ab})$ by minimality of $G$. Since $K^{ab}$ is an abelian, unramified outside of $p$ extension of $\mathbb{Q}$, it is contained in $\mathbb{Q}(\mu_{p^\infty})$ by class field theory. In particular, there is a unique prime $\pi$ of $O_{K^{ab}}$ over $p$. Clearly $H(\pi) = H(\beta) = G(\beta)$ for all primes $\beta$ over $p$, since $\pi = \beta \cap O_{K^{ab}}$. In other words, $G(\beta)$ does not depend on the prime $\beta$ above $p$. Thus we have $H = G(\beta)$ for all $\beta$ over $p$ and therefore all the inertia groups of primes over $p$ act trivially on $G/H$. Consequently $\Gamma$ acts trivially on $G/H$ by Lemma 5. In particular, for any $A \in G$ and any $\tau \in \Gamma$ there is an $h \in H$ such that $A^\tau = hA$. Note that Lemma 4 allows us to assume that $H$ consists of diagonal matrices (after conjugating by a matrix with entries in $\mathbb{Z}$, if necessary). By the following Lemma 7 we get that $A \in GL_n(K^{ab})$ and therefore that $G \subseteq GL_n(K^{ab})$, which contradicts our assumption. $\square$

**Lemma 7** *Suppose that $K$ is a Galois extension of $\mathbb{Q}$ unramified outside a prime $p$. Let $A \in GL_n(O_K)$ has the property that for any $\tau \in Gal(K/\mathbb{Q})$ there is a diagonal matrix $D_\tau = diag(\xi_1(\tau), ..., \xi_n(\tau))$ such that $\xi_i(\tau)$'s are roots of unity of $p$−power order and $A^\tau = D_\tau A$. Then $A \in GL_n(K^{ab})$.*

*Proof:* Let $A = (a_{i,j})$. Let $\beta$ be a prime of $O_K$ over $p$. For every $i$ there exists $k$ such that $a_i = a_{i,k}$ is not in $\beta$. We have $a_{i,j}^\tau = \xi_i a_{i,j}$ for all $j$. In particular, $a_i a_{i,j}^{-1}$ is stable under the action of $Gal(K/\mathbb{Q})$ for all $j$. Thus there are rational numbers $q_{i,j}$ such that $a_{i,j} = a_i q_{i,j}$.

Since $a_i^\tau = \xi_i(\tau) a_i$, there is an integer $n$ such that $(a_i^{p^n})^\tau = a_i^{p^n}$ for all $\tau \in Gal(K/\mathbb{Q})$. Thus $a_i^{p^n} = l_i \in \mathbb{Q}$ and consequently they are integers prime to $p$. We get that $\sqrt[p^n]{l_i} \in K$. Since $K$ is unramified outside $p$ and $l_i$ is prime to $p$, we get that $l_i = \pm t_i^{p^n}$ for some integers $t_i$. This implies that $a_i = t_i \zeta_i$ for some roots of unity $\zeta_i$ in $K$. It follows that $A = diag(\zeta_1, ..., \zeta_n)Q$, where $Q = (t_i q_{i,j}) \in GL_n(\mathbb{Z})$. $\square$

**Proposition 2** *$G$ is an elementary abelian $p$−group.*

*Proof:* It follows from Lemma 6 that $G$ is a $p$−group generated by the subgroups $G(\beta)$ for primes $\beta$ over $p$. Let $W = G \cap GL_n(O_{K^{ab}})$. Every proper, $\Gamma$−stable subgroup of $G$ is contained in $W$ by minimality of $G$. In particular, we have the inclusion $[G, G] \subseteq W$ and therefore $W$ is a normal subgroup of $G$. Since there is a unique prime of $O_{K^{ab}}$ over $p$, the subgroup $H = G(\beta) \cap W$ is independent on the prime $\beta$ above $p$. Thus it is a $\Gamma$−stable, normal subgroup of $G$ and by Lemma 4 we can assume that it consists of diagonal matrices. For a prime $\beta$ over $p$ let $G'(\beta)$ be the image of $G(\beta)$ in $G/H$. Since $[G(\beta), G(\beta_1)] \subseteq G(\beta) \cap G(\beta_1) \cap [G, G] \subseteq H$, we get that $[G'(\beta), G'(\beta_1)] = 1$. In other words, these groups pairwise commute and are commutative. Let $F$ be the Frattini subgroup of $G$. Since $F$ is characteristic, it is $\Gamma$−stable and consequently $F \subseteq W$. Note that $F \cap G(\beta) \subseteq W \cap G(\beta) = H$. Thus the groups $G'(\beta)$ have trivial intersection with the image of $F$ in $G/H$. It follows that they are elementary abelian $p$−groups. Consequently, these subgroups generate elementary abelian $p$−group. On the other hand they generate whole $G/H$, so $G/H$ is an elementary abelian $p$−group generated by $G'(\beta)$, where $\beta$ runs over primes above $p$. It remains to show that $H = 1$.

Suppose that $H$ is not trivial. Thus it is a $p$−subgroup of the diagonal matrices of $GL_n$ and there is no harm to assume that it is in a strongly diagonal form. We can associate to $H$ its normalizer $N(H)$, centralizer $C(H)$ and representation $\rho$, as discussed in Section 2. Clearly $G \subseteq N(H)$. Note that the kernel of $\rho$ is contained in the group of diagonal matrices and the subgroup of diagonal matrices in $G$ is exactly $H$ (such matrices are clearly in the congruence subgroups and have entries in $K^{ab}$). Thus $\rho(G)$ and $G/H$ are isomorphic. By Lemma 3 the image $\rho(G)$ is a finite, $\Gamma$−stable subgroup of $GL_m(O_K)$ for some $m \geq n$, Moreover, the action of $\Gamma$ on $\rho(G)$ coincides with its action on $G/H$ induced from the action on $G$, i.e $G/H$ and $\rho(G)$ are isomorphic as $\Gamma$−modules. Directly from the definition of $\rho$ we get that $\rho(G(\beta)) \subseteq \rho(G)(\beta)$ for any prime $\beta$ of $O_K$. In particular, $\rho(G)$ is generated by the subgroups $\rho(G)(\beta)$ for primes $\beta$ over p. By minimality of our counterexample $G$ we get that the commutator of $\Gamma$ acts trivially on $\rho(G)$. In other words, $\rho(G)$ is contained in $GL_m(K^{ab})$. Since there is a unique prime $\pi$ in $K^{ab}$ over $p$, we have the equality $\rho(G) = \rho(G)(\pi)$. Thus $\rho(G)$ can be conjugated to diagonal matrices over $\mathbb{Z}$ by Lemma 4.

Now we will analyze the action of $\Gamma$ on $H$ and $G/H$

First note that $p$ is odd. In fact, we already proved that $[\Gamma, \Gamma]$ acts trivially on both $H$ and $G/H$. It follows that the commutator $[\Gamma, \Gamma]$ is an abelian $p$−group. In fact, for any $g \in G$ the map $\tau \mapsto g^\tau g^{-1}$ is a homomorphism from $[\Gamma, \Gamma]$ to $H$, and all these homomorphisms separate elements of $[\Gamma, \Gamma]$. If $p = 2$ then the abelianization $\Gamma/[\Gamma, \Gamma]$ is a 2−group ($K^{ab}$ is an abelian, unramified outside 2 extension of $\mathbb{Q}$) and therefore $\Gamma$ itself is a 2−group. Since 2−groups are nilpotent, we get a contradiction by Theorem 1. Thus in the case of $p = 2$ the group $H$ has to be trivial.

Assume now that $p$ is odd. Let $\mu$ be the subgroup of $K^\times$ consisting of roots of unity of $p$−power order. $\Gamma$ acts on $\mu$ via the cyclotomic character $i$. To be more precise, $i(\tau)$ is the unique integer (mod $|\mu|$) such that $\xi^\tau = \xi^{i(\tau)}$ for all $\xi \in \mu$. Since $H$ consists of diagonal matrices, $\Gamma$ acts on $H$ via $i$, i.e. $h^\tau = h^{i(\tau)}$ for all $h \in H$. The same is true for $\rho(G)$, since it is conjugated to diagonal matrices over $\mathbb{Z}$. Thus, $b^\tau = b^{i(\tau)}$ for all $b \in G/H = \rho(G)$ and $\tau \in \Gamma$. As a consequence we get that for any $g \in G$ the subgroup of $G$ generated by $g$ and $H$ is $\Gamma$−stable. Let $g \in G$ be any element which is not fixed by $[\Gamma, \Gamma]$. The subgroup $< H, g >$ of $G$ generated by $g$ and $H$ is $\Gamma$−stable and the commutator of $G$ acts nontrivially on it. By minimality of $G$ we get that $G = < H, g >$. In particular, $G/H$ is a cyclic group of order $p$. We need to consider two cases:

*case 1:* $H$ is central in $GL_n$.

Thus $H$ is cyclic, $G$ is abelian and either $G$ is cyclic or $G = H \times \mathbb{Z}/p$. The former case implies that $[\Gamma, \Gamma]$ acts trivially on $G$, which is impossible. Thus $G = H \times F$, where $F$ is cyclic of order $p$. We can assume that $g$ is a generator of $F$. For any $\tau \in \Gamma$ we have $\tau(g) = u_\tau g^{i(\tau)}$ for some $u_\tau \in H$. Clearly $u_\tau$ is of order $p$. By minimality of $G$ it follows that $H$ is generated by the $u_\tau$'s. In particular, $H$ has exponent $p$. Since $H$ is cyclic, it is of order $p$. Let $S$ be the subgroup of $\Gamma$ which acts trivially on both $H$ and $G/H$. Note that $\tau^{p-1} \in S$ for any $\tau \in \Gamma$, because both $H$ and $G/H$ are cyclic of order $p$. Since the commutator of $\Gamma$ is contained in $S$, the group $\Gamma/S$ is cyclic of order dividing $p - 1$. The order has to be exactly $p - 1$ by looking at the action on $p$−th roots of unity. On the other hand, the map $s^{-1} \mapsto g^s g^{-1}$ is an injective homomorphism of $S$ into $H$. Thus $S$ is of order $p$. Putting these together we conclude that $S = [\Gamma, \Gamma]$ is cyclic of order $p$ and $\Gamma/S = \text{Gal}\mathbb{Q}(\zeta)/\mathbb{Q}$ is cyclic of order $p - 1$, where $\zeta$ is a primitive $p$−th root of 1.

Now let $\tau \in \Gamma$ maps to a generator of $\Gamma/S$. Thus $\tau$ acts on $p - th$ roots of unity by raising to a power $t$, where $t$ is a primitive root mod $p$. We have $g^\tau = ug^t$ and $g^{\tau^i} = u^{it^{i-1}} g^{t^i}$ for some $u \in H$ and any $i \in \mathbb{Z}$. If $u \neq 1$ then $\tau$ has order $p(p-1)$ and $\Gamma$ is cyclic- a contradiction. Thus $u = 1$ for any $\tau$ mapping to a generator of $\Gamma/S$. But such elements generate $\Gamma$, so $\Gamma$ preserves $F$. But then the commutator of $\Gamma$ acts trivially on $F$ and $H$, hence on $G$, again a contradiction.

*case 2:* $H$ is not central in $GL_n$.

Then the centralizer $C(H)$ equals $GL_{a_1} \times ... \times GL_{a_s}$, where $a_i < n$ for all $i$. Since $Z = G \cap C(H)$ is $\Gamma$–stable, we get that $[\Gamma, \Gamma]$ acts trivially on $Z$ by minimality of $n$. To be more precise, the projection $Z_i$ of $Z$ into $GL_{a_i}(O_K)$ is $\Gamma$–stable and either $Z_i$ has order smaller than the order of $Z$ or $Z_i$ is isomorphic to $Z$ but sits in matrices of smaller dimension. In both cases $[\Gamma, \Gamma]$ acts trivially on $Z_i$ by minimality of $G$ and $n$. Since this is true for all $i$, it follows that $[\Gamma, \Gamma]$ acts trivially on $Z$. Clearly $G \subseteq N(H)$. Thus we can write $g = cw$ with $c \in C(H)$ and $w$ a permutation matrix. Since $\Gamma$ acts on $G/H$ via the cyclotomic character, there is $\tau \in \Gamma$ such that $g^\tau = ug^2$ for some $u \in H$. On the other hand, we have $g^\tau = c^\tau w$ and consequently $ucwcw = c^\tau w$, so $w = (uc)^{-1}c^\tau c^{-1} \in C(H)$. Thus $g \in Z$ and therefore $G = Z$ has trivial $[\Gamma, \Gamma]$–action contrary to our assumption.

Both cases have led to a contradiction and therefore our assumption that $H$ is not trivial has to be wrong. □

We summarize our investigation in the following theorem:

**Theorem 2** *Suppose that Conjecture 4 is true in the following very special situation:*

*— $K$ is a Galois extension of $\mathbb{Q}$ unramified at all finite primes $q \neq p$;*

*— $G$ is an elementary abelian $p$–group contained in $GL_n(O_K)$ and stable under the action of $GalK/\mathbb{Q}$;*

*— $G(\beta) \cap GL_n(K^{ab}) = 1$ for all primes $\beta$ over $p$ and these subgroups generate $G$.*

*Then it is true in general.*

**Remark.** Recently Kitaoka proved ([5]) that Conjecture 4 is true for $n = 2$. The main step in the proof is the reduction to abelian groups $G$. Note however that our reduction can not be used for this purpose. The point is that we used our "dirty trick" of raising the dimension. It would be nice to have a proof of Theorem 2 which does not use any such methods. Unfortunately we were unable to produce any satisfactory argument in general. However for totally real fields one can avoid raising the dimension for groups $G$ of odd order. The key point is that we have Corollary 1 which gives triviality of the congruence subgroups over $\mathbb{Z}$. It is worth pointing out that in the totally real case (or more generally, when complex conjugation is central) it is known that Conjecture 4 holds for $n < 43$ (see [1]).